\documentclass[11pt]{article}
\usepackage{graphicx}
\usepackage{amssymb, amsmath}

\newtheorem{theorem}{Theorem}[section]

\newtheorem{definition}[theorem]{Definition}
\newtheorem{lemma}[theorem]{Lemma}
\newtheorem{Remark}[theorem]{Remark}
\newtheorem{Proposition}[theorem]{Proposition}
\title{\textbf{Codes defined by forms of degree 2 on hermitian surfaces and S\o rensen's conjecture}}          
\date{}
\author{\textbf{Fr\'ed\'eric A. B. Edoukou}  \\ \\ CNRS, Institut de Math\'ematiques de Luminy  \\ Luminy case 907  - 13288  Marseille Cedex 9 - France \\  E.mail : edoukou@iml.univ-mrs.fr }
\begin{document}
\maketitle

\begin{flushleft}
\textbf{Abstract}
\end{flushleft}
We study the functional codes $C_h(X)$ defined by G. Lachaud in $\lbrack 10 \rbrack$ where $X \subset {\mathbb{P}}^N$ is an algebraic projective variety of degree $d$ and dimension $m$.\\
When $X$ is a hermitian surface in $PG(3,q)$, S\o rensen in \lbrack 15\rbrack, has conjectured for $h\le t$ (where $q=t^2$) the following result :
$$\# X_{ Z(f)}(\mathbb{F}_{q} ) \le h( t^{3}+ t^{2}-t)+t+1$$
which should give the exact value of the minimum distance of the functional code $C_h(X)$.\\
In this paper we  resolve the conjecture of S\o rensen in the case of quadrics  (i.e. $h=2$), we show the geometrical structure of the minimum weight codewords and their number; we also estimate the second weight  and the geometrical structure of the codewords reaching this second weight.\\\\
\noindent \textbf{Keywords:} functional codes, hermitian surface, hermitian curve, quadric, S\o rensen's conjecture, weight.\\\\
\noindent \textbf{Mathematics Subject Classification:} 05B25, 11T71, 14J29
\section{Introduction}
After the works of Goppa, in coding theory, the construction of error-correcting codes from algebraic curves is now classical.  From the works of renowned scholars in algebraic geometry, such as Manin, Vladut $\lbrack 11 \rbrack $, we are able to consider codes built from higher dimensional algebraic varieties. \\
Some of such codes have already been studied. Among others, projective Reed-Muller codes have been studied by G. Lachaud $\lbrack 9 \rbrack $, A. B. S\o rensen $\lbrack 14 \rbrack $ and Y. Aubry  $\lbrack 1 \rbrack $;  codes on hermitian varieties have also been studied by Chakravarti $\lbrack 4 \rbrack $, and Hirschfeld,  Vladut, Tsfasman  $\lbrack 8 \rbrack $.\\ 
In $\lbrack 10 \rbrack $  G. Lachaud has also study the functional codes on a projective variety $X$ of $\mathbb{P}^{N}$ of degree $d$ and dimension $m$. Theses codes have similarities with the well known projective Reed-Muller codes defined over a projective variety $X$ of $\mathbb{P}^{N}$.\\
In this paper, what we have in mind is to study the functional codes $C_{h}(X)$ introduced by G. Lachaud in $\lbrack 10 \rbrack$  in the particular case where $X$ is the non-degenerate hermitian surface $X : x_{0}^{t+1}+x_{1}^{t+1}+x_{2}^{t+1}+x_{3}^{t+1}=0$ in $PG(3, q) $ and $q=t^2$. Indeed, the hermitian surfaces are interesting since they are maximal with respect to their number of rational points and they have also a beautiful geometry. \\
Even in the case of quadric surfaces,  Y. Aubry $\lbrack 1 \rbrack $ was not able to determine precisely the minimum distance of the functional codes $C_2(X) $, when the quadrics are non-degenerate (elliptic or hyperbolic) and 1-degenerate (cones). \\
S\o rensen $\lbrack 15 \rbrack $ has also tried to study the functional codes $C_h(X)$, on the hermitian surface $X$; but he was not able to find the minimum distance, even in the case $h=2$.  His work was motivated by the paper $\lbrack 3 \rbrack $ (of Chakravarti) which undertake the study of the code $C_2(X)$ over  $\mathbb{F}_4$ by a complete computer search.  Thus, he  conjectured a result, which  will be the main purpose of our work. \\
The paper has been organized as follows. First of all, we recall somme notations and definitions of functional codes. Secondly, we state the conjecture of S\o rensen: $$\mathrm{for}\  h\le t,\quad \# X_{ Z(f)}(\mathbb{F}_{q} ) \le h( t^{3}+ t^{2}-t)+t+1.$$
Thirdly, by using the projective classification of quadrics in  $\lbrack 7 \rbrack $, the properties of the hermitian surface $X$ and the projective classification of the hermitian plane curves in $\lbrack 7 \rbrack $,  we give a proof of the conjecture of S\o rensen in the case $h=2$. Finally we conclude our work by expressing more clearly the exact parameters of the functional codes $C_2(X)$, the geometrical structure of the minimum weight codewords and their number, the second weight and also the geometrical structure of codewords reaching this second weight. 
\section{Notations}
We denote by $\mathbb{F}_{q} $ the field with  $q$ elements, where $q=p^{a} $ is a power of the prime number $p$. Let  $V=\mathbb{A}^{N+1}$ be the affine space of dimension $N+1$ over $\mathbb{F}_{q}$. We denote also by $\mathbb{P}=\mathbb{P}(V)=PG(N,q)={\mathbb{P}}^{N}$ the corrresponding projective space of dimension equal to $N$.\\
We denote also by $P_{h}(V,\mathbb{F}_{q})$ the vector space of forms (that is, of homogeneous  polynomials) of degree $h$ in $V$ with coefficients in $\mathbb{F}_{q}$. If $f \in P_{h}(V,\mathbb{F}_{q})$, we denote by $Z(f)$, the set of zeros of $f$ in the projective space. Let $X\subset\mathbb{P}$ be a subvariety of $\mathbb{P}$, $\vert X\vert$ denotes the number of its rational points over $\mathbb{F}_{q}$; we call $X\cap Z(f)$ a degree $h$ section of $X$, and we denote by ${X_{Z(f)}}(  \mathbb{F}_{q}) $, the the set of rational points over $\mathbb{F}_q$ of the algebraic set $X\cap Z(f)$.
 \section{Functional codes defined by forms on projective varieties}
The subject of this section is to recall the construction of the functional codes defined by forms on projective varieties as it has been done by G. Lachaud in  $\lbrack 10 \rbrack $. 
 \begin{definition}
Let $X$ be a finite set, $X=(P_{1},\dots,P_{n})$.\\
Let $\mathcal{F}(X, \mathbb{F}_{q})$ the space of all maps from $X$ to $\mathbb{F}_{q}$.      $\mathcal{F}(X, \mathbb{F}_{q})$ is a vector space; let $\widetilde{F}\subset \mathcal{F}(X, \mathbb{F}_{q})$ a subspace.\\
Let  $c$ be the map defined  by $$
 c:\mathcal{F}(X, \mathbb{F}_{q})
  \longrightarrow
  \mathbb{F}_{q}^{n}$$
  $$\qquad \qquad \qquad \qquad \qquad \qquad \qquad  f\longmapsto c(f)=(f(P_{1}),\ldots,f(P_{n}))
  $$
We call functional code defined by $\widetilde{F}$ and $X$, and we denote it by $C(X,\widetilde{F})$, the image of the map $c$ restricted to $\widetilde{F}$.
 $$
 \qquad \quad \  c_{\vert \widetilde{F}}:\mathrm{\widetilde{F}}
  \longrightarrow
 \mathbb{F}_{q}^{n}$$
  $$\qquad \qquad \qquad \qquad \qquad \qquad \qquad  \qquad f\longmapsto c_{\vert \widetilde{F}}(f)=(f(P_{1}),\ldots,f(P_{n}))
  $$
$$C(X,\widetilde{F})=\mathrm{Im} c_{\vert \widetilde{F}} $$
The functional code we have defined has the following parameters
$$\mathrm{lengh} C(X,\widetilde{F})=n,\qquad  \dim C(X,\widetilde{F})=\dim \widetilde{F}-\dim \ker C_{\vert \widetilde{F}}$$

\end{definition}
\begin{definition}
Let $c(f)$ be a codeword. We call coweight of the word $c(f)$, and we denote $cw(f)$, the number of zeros digits in  $c(f)$.\\
$$ cw(f)= \# \{P \in X\ \  \vert \quad f(P)=0 \} \qquad (3.2.1)$$
The weight of the word $c(f)$ which we denote by $w(c(f))$ is 
 $$w(c(f))=n-cw(f) $$
 The minimal distance of the functional code $C(X,\widetilde{F})$ is 
  $$\mathrm{dist}C(X,\widetilde{F})=n-\underset{\mathrm{f}\in \mathrm{\widetilde{F}}}{\max}\ cw(f) \quad (3.2.2)$$
\end{definition}

\begin{definition}
We denote by $( x_{0}: \ldots: x_{N})$ the homogeneous coordinates in the projective space 
$\mathbb{P}^{N} $.\\\
Let  $W_{i}$ be the set of points with homogeneous coordinates $(x_{0}: \ldots: x_{N})  \in \mathbb{P}^{N}$ such that $x_{0}=x_{1}=\ldots x_{i-1}=0$ and $x_{i}\ne 0$.\\
The family $ \{ W_{i} \}_{ 0 \leq i \leq N } $ is clearly a partition of  $\mathbb{P}^{N} $.\\
Let  $\nu$ be a map defined by: $$\quad  \nu : \mathbb{P}^{N} \longrightarrow V-\{0\}$$
$$\qquad \qquad \quad
( x_{0}:\ldots: x_{N}) 
 \longmapsto 
 (0,\ldots,0,1,x_{i+1},\ldots,x_{N})
 $$
for  $( x_{0}:\ldots: x_{N}) \in W_{i}$.\\\\ 
If $f \in P_{h}(V,\mathbb{F}_{q})$ and $Q \in \mathbb{P}^{N}$, then
$$f(\nu Q)=\frac{f(x_{0},\ldots,x_{N})}{x_{i}^{h}} \qquad  \mathrm{   with  \qquad Q\in    W_{i}}$$
\end{definition}
\begin{definition}
Let $X \subset \mathbb{P}^{N}$ be an algebraic set defined over $ \mathbb{F}_{q}  $ and $h\leq q$.
The code defined by forms of degree $h$ on $X$ is the functional code  $C_{h}(X)=C(X,\widetilde{F})$  obtained by taking $\widetilde{F}$ as the space of restrictions to $X$  of forms in $P_{h}(V,\mathbb{F}_{q})$.

\end{definition}
We have  $$\mathrm{length}\ C_{h}(X) = \# X(\mathbb{F}_{q} ) \quad (3.4.1)$$  \\
When the map $c$ is injective, we get  $$ \dim C_{h}(X)= {N+h\choose h} \quad (3.4.2)$$
From $(3.2.1)$, we deduce that  
$$cw(f)=\# X_{ Z(f)}(\mathbb{F}_{q} ) \quad (3.4.3) $$
From $(3.4.3)$, we can conclude that calculating the weights of the code $C_{h}(X)$, comes down to the computation of the number of points of degree $h$ sections of $X$.
\section{Codes defined by forms on hermitian varieties in $PG(3,q)$}
We now restrict ourselves to functional codes defined in particular algebraic sets in $PG(3,q)$, where  $q=t^2$ and $t$ is  a prime power. Thus, we consider the hermitian surface $X$ defined by  :
$$ X : x_{0}^{t+1}+ x_{1}^{t+1}+x_{2}^{t+1}+ x_{3}^{t+1}=0  $$
We know from $(3.4.1)$ that the length of this code is   $\mathrm{length}\ C_{h}(X) = \# X(\mathbb{F}_{q})$. In \lbrack 2\rbrack  (p.1179)  Bose and  Chakravarti showed that  $$ \# X(\mathbb{F}_{q} )=(t^2+1)(t^3+1)\quad (4.1 ) .$$
If $h\le t$ the map c is injective and from $(3.4.2)$, we can say that $$ \dim C_{h}(X)= {3+h\choose h} \quad (4.2).$$
For an estimate on the minimum distance of this code, as usually, a more careful analysis is necessary.
\subsection{Codes defined by linear forms}
In $\lbrack 2 \rbrack$  we have the following result.
\begin{theorem}
Let $H \subset \mathbb{P}^3 $ be a plane
\begin{equation*}
 \# X_{H}(\mathbb{F}_{q})=
         \begin{cases}
            t^3+1 &    \text{\  \ if\  H\  is\  not\  tangent \ to \ X,} \\
            t^3+t^2+1 &  \text{ \ if \ H \ is \ tangent \ to \ X}.
            \end{cases} 
\end{equation*}                         
\end{theorem}

If $H$ is tangent to $X$, $X\cap H$ is a singular hermitian curve of rank $2$ in $PG(2, t^2)$. More precisely, $H\cap X$ is a set of t+1 lines passing through a common point. Thus we get $\# X_{H}(\mathbb{F}_{q})= (t+1)t^2+1=t^3+t^2+1.$ \\
If $H$ is not tangent to $X$, $X\cap H$ is a non-singular hermitian curve in $PG(2, t^2)$. Thus, we get $\# X_{H}(\mathbb{F}_{q})= t^3+1$.
\subsection{Codes defined by forms of degree two (quadrics)}
In the case  $h=2$ and $t=2$,  the minimum distance of this code has been found by Paul P. Spurr in 1986 in his master's project by using a complete computer search;  thus a $ \lbrack 45,10, 22 \rbrack$-code over $\mathbb{F}_{4} $ has been found.\\
In  1991, Anders B. S\o rensen in $ \lbrack 15 \rbrack$, in his Ph. D. Thesis, without any computer programs, only with the geometric properties of quadrics and those of the non-singular cubic in $\lbrack 5 \rbrack$, found the minimum distance of this code defined over $\mathbb{F}_{4} $.\\
S\o rensen has tried to generalize the problem and conjectured that for  $h \le t$ (and $t$ any prime power) we get 

 $$\# X_{ Z(f)}(\mathbb{F}_{q} ) \le h( t^{3}+ t^{2}-t)+t+1  $$

G. Lachaud has also proved in $\lbrack 10 \rbrack$, with the help of proposition 2.3 $\lbrack 10 \rbrack $ that   $$ \# X_{ Z(f)}(\mathbb{F}_{q} ) \le h( t^{3}+ t^{2}+t +1)$$
Unfortunately his bound is larger than that of S\o rensen.
\begin{Remark}
From theorem 4.1 the conjecture is true in the case of codes defined by linear forms (h=1). In that case, the code $C_1(X)$ (noted  $C(X)$) is a two-weight code of dimension $4$ whose weights are $t^5$ and $t^5+t^2 $.
\end{Remark}
We will now resolve the problem of Sorensen, in the case $h=2$, and $t$ any prime power.
\section{Resolution of S\o rensen's conjecture in the case $h=2$ and $t$ any prime power}

We need to recall the notion of the rank of a quadric $\mathcal{Q}=Z(f)$. 
\begin{definition}
Let $T$ be an inversible linear transformation defined over $\mathbb{P}^{N}$.
Let $i(f)$ be the number of indeterminates appearing explicitly in $F$. The rank $r(\mathcal{Q})$ of the quadric $\mathcal{Q}$ is defined by $$r(Q)=\underset{T}{\min}\ i(f\circ T)$$
where T ranges over all inversible linear transformations over $\mathbb{F}_q$.
A quadric $\mathcal{Q}$ is said to be degenerate if  $r(\mathcal{Q})< N+1$. Otherwise the quadric is non-degenerate.
\end{definition}
We have a particular property of quadric in \lbrack 6\rbrack   \   p.115.
\begin{Proposition}
  A quadric $\mathcal{Q}$ is degenerate if and only it is singular. 
\end{Proposition}
We recall also the classification of quadrics in  $PG(3,q)$ under  $PGL(4,q)$ given by  J.W.P. Hirschfeld   $\lbrack 7 \rbrack$ in table 1; thus we will consider two cases depending on the degeneration of the quadrics.\\
In the whole  paragraph we will note 
$$ s(t)=2(t^3+t^2-t)+t+1  \qquad  s_{2}(t)=2t^3+t^2+1$$

\unitlength=1cm
\hspace{-13 mm}
\begin{picture}(12,7.5)
\put(2,0){\line(1,0){10}}
\put(2,1){\line(1,0){10}}
\put(2,2){\line(1,0){10}}
\put(2,3){\line(1,0){10}}
\put(2,4){\line(1,0){10}}
\put(2,5){\line(1,0){10}}
\put(2,6){\line(1,0){10}}
\put(2,7){\line(1,0){10}}

\put(2,0){\line(0,1){7}}
\put(5,0){\line(0,1){7}}
\put(9,0){\line(0,1){7}}
\put(12,0){\line(0,1){7}}
\put(3.2,6.3){Rank}

\put(5.3,6.3){Description}
\put(9.3,6.3){$ \vert \mathcal{Q} \vert $}
\put(3.5,5.3){1}
\put(3.5,4.3){2}
\put(3.5,3.3){2}
\put(3.5,2.3){3}
\put(3.5,1.3){4}
\put(3.5,0.3){4}

\put(9.3,5.3){$q^2+q+1$}
\put(9.3,4.3){$2q^2+q+1$}
\put(9.3,3.3){$q+1$}
\put(9.3,2.3){$q^2+q+1$}
\put(9.3,1.3){$(q+1)^2$}
\put(9.3,0.3){$q^2+1$}

\put(5.3,5.3){repeated plane} 
\put(5.3,4.3){pair of distinct planes }
\put(5.3,3.3){line }
\put(5.3,2.3){(quadric) cone}
\put(5.3,1.3){hyperbolic quadric}
\put(5.3,0.3){elliptic quadric}

\put(4.5,7.3){ Table.1  Quadrics in PG(3,q) }

\end{picture}
\subsection{$\mathcal {Q } $ is a degenerate quadric}

Degenerate quadrics in $PG(3,q)$ are quadrics with rank stricly less than $4$.\\ 
So we need to distinguish four cases.
\subsubsection{Rank $\mathcal {Q } =1$ }
In this case $\mathcal {Q } $ is a repeated plane.\\
If $\mathcal {Q }$ is tangent to $X$, from theorem 4.1, $\mathcal {Q } \cap X$ is a singular hermitian  curve of rang $2$ in $PG(2,t^2)$, and $\vert \mathcal{Q }\cap X\vert =t^3+t^2+1< s_2(t)$.\\
If $\mathcal {Q } $ is not tangent to $X$,  from $ \lbrack 2 \rbrack$, $\mathcal {Q } \cap X$ is a non-singular hermitian curve  in  $PG(2,t^2)$, and $\vert\mathcal{Q }\cap X\vert=t^3+1 < s_2(t).$

\subsubsection{Rank $\mathcal {Q } =2$ and $\mathcal {Q }$  is a line}
In this case we have $\vert \mathcal {Q } \vert=q+1=t^2+1$, so $\vert \mathcal{Q }\cap X\vert \le t^2+1< s_2(t)$.

\subsubsection{Rank $\mathcal {Q } =2$ and $\mathcal {Q }$  is a pair of distinct planes}
\paragraph{(i) each plane is not tangent to $X$}. \\
If we write $\mathcal {Q } =  \mathcal {P} _{1} \cup   \mathcal {P} _{2}$, then $\mathcal {P }_{1} \cap X $ and $\mathcal {P }_{2} \cap X $ are non-singular hermitian curves in $PG(2,t^2)$, such that $\vert  \mathcal {P }_{1} \cap X \vert  = \vert \mathcal {P }_{2} \cap X \vert =t^3+1$. Hence $\vert \mathcal{Q }\cap X \vert \le 2(t^3+1) < s_2(t)$.

\paragraph{(ii) one plane is tangent and the second is not tangent to $X$}.\\
Let us note $\mathcal {Q } =  \mathcal {P} _{1} \cup   \mathcal {P} _{2}$, with $\mathcal {P} _{1}$ tangent to $X$ and $\mathcal {P} _{2}$ not tangent to $X$. The singular hermitian curve of rang $2$, $\mathcal {P }_{1} \cap X $ is in fact a set of $(t+1)$ lines passing through a point $P$, see theorem 4.1.\\
Let $l=\mathcal {P} _{1} \cap   \mathcal {P} _{2}  $, then $l$ is not contained in $X$ from theorem 9.1 $\lbrack 3\rbrack$ p.1176. Since $l \subset \mathcal {P} _{1}$, we have $ l \cap X=(l \cap \mathcal {P} _{1} )\cap X =l \cap (\mathcal {P} _{1} \cap X)$, so that $l \cap X $ is a single point, or a set of $t+1$ points since two distinct lines in the projective plane have exactly one common point.  
We know that 
 $$\vert \mathcal{Q }\cap X \vert = \vert  \mathcal {P }_{1} \cap X \vert  +\vert \mathcal {P }_{2} \cap X \vert - \vert \mathcal {P} _{1} \cap   \mathcal {P} _{2}  \cap X \vert   \quad (5.1.3.\mathrm{ii })$$ 
 Since $\vert  \mathcal {P }_{1} \cap X \vert= t^3+t^2+1$ and  $\vert \mathcal {P }_{2} \cap X \vert=t^3+1$, we get  $\vert \mathcal{Q }\cap X \vert=2t^3+t^2+1=s_2(t)$ in the case $l\cap X$ is a single point, or $\vert \mathcal{Q }\cap X \vert=2t^3+t^2-t+1<s_2(t)$ in the case $l\cap X$ is a set of $t+1$ points.

\paragraph{(iii) each plane is tangent to $X$}.\\
Let us note $\mathcal {Q }= \mathcal {P} _{1} \cup \mathcal {P} _{2}$, with $\mathcal {P} _{1}$, and $\mathcal {P} _{2}$  tangent to $X$.  We know that  $\mathcal {P }_{1} \cap X $ and $\mathcal {P }_{2} \cap X $ as singular hermitian curves in $PG(2,t^2)$  are each one a set of  $t+1$ lines passing through respectively a point $P_{1}$ and a point $P_{2}$. Let $l =  \mathcal {P} _{1} \cap   \mathcal {P} _{2}  $, then $l \cap X$ is a set of $ t+1$ points or a line (in the case $l$ is contained in $X$).  Since  $ \vert  \mathcal {P }_{1} \cap X \vert  = \vert \mathcal {P }_{2} \cap X \vert =(t+1)t^2+1$, from the relation $(5.1.3.\mathrm{ii})$ we conclude that  either $\vert \mathcal{Q }\cap X \vert= s(t)$ when $l \cap X$ is a set of $ t+1$ points, or $\vert \mathcal{Q }\cap X \vert=2t^3+t^2+1= s_2(t)$ when $l$ is contained in $X$.
\subsubsection{Rank $\mathcal {Q } =3$ ($\mathcal {Q }$  is a cone)}
Here, $\vert \mathcal {Q } \vert=q^2+q+1=q(q+1)+1$ and 
$\mathcal {Q }$ consists of the points on $q+1$ lines, passing through a vertex $S$. Thus we can consider two cases.

\paragraph{(i)  no line in the cone is contained $X$}.\\
Let us note $\mathcal {Q } =\bigcup _{i=1}^{q+1}\mathcal{L}_{i} $ where  $\mathcal{L}_{i} $ is a line
$$  \mathcal {Q } \cap X = \bigcup _{i=1}^{q+1}(X \cap \mathcal{L}_{i} ) \quad \mathrm{(5.1.4.i.1)}$$
We have also  $ \vert X\cap \mathcal{L}_{i}\vert  \le t+1$, so from the relation $\mathrm{(5.1.4.i.1)}$ $\vert \mathcal{Q}\cap X \vert \le (t+1)(q+1)$. Hence $\vert \mathcal{Q }\cap X \vert \le t^3+t^2+t+1 <  s_2(t)$.

\paragraph{(ii)  at least one  line of the  cone is contained in  $X$}.\\
In $PG(3, q)$ the number of lines passing through a point is  $q^2+q+1=t^4+t^2+1$.\\
From the study of the hermitian surface $X$ in $\lbrack 2\rbrack$, we have the following result.
\begin{Proposition}
If C is any point in $X$, there pass exactly $t+1$ generators (i.e. lines lying in $X$) which constitute the intersection with $X$ of the tangent plane at $C$. And through $C$, there pass $t^2+1$ lines lying in the tangent plane at $C$, among them $t+1$ are generators. The remaining $t^2-t$ lines through $C$, which lie in the tangent plane meet $X$ only in the single point $C$ (there are called tangents to $X$). Each of the $t^4$ lines passing through $C$ and not contained in the tangent plane at $C$ intersect $X$ in exactly $t+1$ points, one of which is $C$ (such lines are called secants to $X$).
\end{Proposition}
From this proposition, on the hermitian surface $X$, we can assert first of all that the maximum number of generators the cone $\mathcal {Q } $ can contained is $t+1$. Secondly, if we suppose that the remaining $t^2-t$ lines  of the cone passing through the vertex $S$ are all secants to $X$, we get 
$$\vert \mathcal{Q }\cap X \vert \le  (t+1)t^2+1+(t^2-t)t.$$
And therefore, we have  $\vert \mathcal{Q }\cap X \vert \le 2t^3+1< s_2(t).$ 

\subsection{ $\mathcal {Q }$  is a non-degenerate quadric }
Non-degenerate quadrics in $PG(3,q)$ are quadrics with rank equal to $4$. Thus, from table 1, p. 7 there are two different  types of such quadrics : hyperbolic quadrics and elliptic quadrics.

\subsubsection{$\mathcal {Q }$ is a hyperbolic quadric }
When  $\mathcal{Q}$ is a hyperbolic quadric, we have $ \vert \mathcal {Q } \vert =  (q+1)^2$. \\
Before giving  the proof of the conjecture of S\o rensen when $  \mathcal {Q }$ is a hyperbolic quadric,  we need to recall some properties of lines, hermitian surfaces and hyperbolic quadrics in $PG(3,q)$.
\begin {lemma} $\lbrack 7\rbrack$
In PG(3,q), the number of lines meeting three skew lines is $q+1$.
\end{lemma}
For a proof of this lemma see $\lbrack 7\rbrack$ p.3  .
\begin{definition}
The set of transversals of three skew lines is called a regulus. So that a regulus consists of $q+1$ skew lines, and if $l_1, l_2$ and $l_3$ are any three of them, it is denoted by $\mathcal{R}(l_1, l_2 , l_3) .$
\end{definition}
We have the following theorem, for the proof see $\lbrack 7\rbrack$ p.23  .
\begin{theorem}
In PG(3,q), $q>2$, let $\{ l_1, l_2 , l_3, l_4 \} $ and $\{ l_1^{\prime}, l_2^{\prime} , l_3^{\prime}, l_4^{\prime} \}$ be two sets of four lines such that any two of the same set are skew and such that $15$ of the $16$ pairs $\{ l_i , l_j^{\prime} \} $meet in a point. Then the last pair also does.
\end{theorem} 
 The above theorem says that the lines meeting $l_1, l_2$ and $l_3$ meet all the lines of $\mathcal{R}(l_1, l_2 , l_3) = \mathcal{R}$ and form a regulus, called the complementary regulus of  $\mathcal{R}$ .\\
From these results we have a good description of the hyperbolic quadric in $PG(3,q)$. 
\begin{theorem}
The hyperbolic quadric  $\mathcal{Q}$ consists of $(q+1)^2$ points, which are all on a pair of complementary reguli. The two reguli are the two systems of generators of $\mathcal{Q}$.  
A hyperbolic quadric whose complementary reguli are $\mathcal{R}$ and  $\mathcal{R}^{\prime} $ is denoted by  $\mathcal{H}( \mathcal{R} , \mathcal{R}^{\prime} )$.
\end{theorem}
In $\lbrack 7\rbrack$ p.123, we also have an important property on the intersection of a hyperbolic quadric and the hermitian surface $X$, in the case the hyperbolic quadric contained at least three lines of a regulus.
\begin{theorem}
If $\mathcal{Q}= \mathcal{H}( \mathcal{R} , \mathcal{R}^{\prime} )$ has three skew lines on $X$, then $X \cap \mathcal{Q}$ consists of $2(t+1)$ lines of  $\mathcal{R}_{0} \cup  \mathcal{R}^{\prime} _{0} $ where $\mathcal{R}_{0} \subset  \mathcal{R}$, $\mathcal{R}^{\prime} _{0}  \subset  \mathcal{R}^{\prime}$ and  $\vert \mathcal{R}^{\prime} _{0}  \vert = \vert \mathcal{R}_{0} \vert=t+1 $.
\end{theorem}
\begin{Remark}
Theorem 5.8 says that, if a hyperbolic quadric contains three skew lines on the non-degenerate hermition surface $X$, then it contains exactly $2(t+1)$ lines of the surface $X$, and $t+1$ lines in each of the two reguli.
\end{Remark}
So that, if  $X \cap \mathcal{Q}= \mathcal{C}$, then $\mathcal{C}= \mathcal{C}_1 \cup \mathcal{C}_2 $ where $\mathcal{C}_1$ and $\mathcal{C}_2$ are respectively the set of simple points, and double points, $\mathcal{C}_2=\{ l \cap l^{\prime} \vert l\in \mathcal{R}_{0},\  l^{\prime}\in \mathcal{R}^{\prime} _{0}  \}.$ We get $\vert \mathcal{C}_1\vert= 2(t+1)(t^2+1-(t+1))$ and $\vert \mathcal{C}_2\vert= (t+1)^2$. \\
Hence $\vert \mathcal{Q }\cap X \vert=\vert \mathcal{C}\vert =2t^3+t^2+1.$\\ 

We can consider only one regulus instead of working in the whole hyperbolic quadric.
Thus, we are now in position to say that the apparently difficult problem on the intersection of an hyperbolic quadric and the hermitian surface is reduced to the four following simple cases.
\paragraph{(i)A regulus of the hyperbolic quadric $\mathcal{H}( \mathcal{R} , \mathcal{R}^{\prime} )$ contains at least three skew lines on the surface $X$ }
In this case, from the above remark we have $\vert \mathcal{Q }\cap X \vert= s_2(t)$. 
\paragraph{(ii)A regulus of the hyperbolic quadric $\mathcal{H}( \mathcal{R} , \mathcal{R}^{\prime} )$ contains exactly two skew lines on the surface $X$ }
We know that the hyperbolic quadric $ \mathcal{Q}$ is generated by $q+1=t^2+1$  lines of a regulus  $\mathcal{R}$. The remaining $t^2-1$ lines of this regulus $\mathcal{R}$ which are not contained in the hermitian surface $X$, each of them, meets the hermitian surface $X$, in at most $t+1$ points. Therefore  we have $\vert \mathcal{Q }\cap X \vert \le 2(t^2+1)+(t^2-1)(t+1)=t^3+3t^2-t+1 < s_2(t)$.
\paragraph{(iii)A regulus of the hyperbolic quadric $\mathcal{H}( \mathcal{R} , \mathcal{R}^{\prime} )$ contains exactly one  line on the surface $X$ }
We will suppose now that $\mathcal{R}$ is the regulus of  the hyperbolic quadric which contains the only line of intersection. In this way, there are $t^2$ lines of $\mathcal{R}$ which are not contained in $X$. Each of them, meeting the hermitian surface in at most $t+1$ points,  therefore  $\vert \mathcal{Q }\cap X \vert \le t^2+1+t^2(t+1)=t^3+2t^2+1 < s_2(t)$. 
\paragraph{(iv) Each regulus of the hyperbolic quadric $\mathcal{H}( \mathcal{R} , \mathcal{R}^{\prime} )$ does not contain any line on the surface $X$ }
The hyperbolic quadric generated by one regulus for instance $\mathcal{R}$, and each line of the regulus meeting the hermitian surface  $X$ in at most $t+1$ points, implies that $\vert \mathcal{Q }\cap X \vert \le (t^2+1)(t+1)=t^3+t^2+t+1 < s_2(t)$, as in the case 5.1.4.i). 
\subsubsection{$\mathcal {Q }$ is an elliptic quadric }
When $\mathcal {Q }$ is an elliptic quadric, we have $\vert \mathcal {Q } \vert =q^2+1=t^4+1$ and no line is contained in $\mathcal {Q }$ (that is the characterization of an elliptic quadric in $PG(3, q)$). \\
We need to distinguish two cases.
\paragraph{(i) $X \cap \mathcal {Q } = \emptyset $}.\\
In this case the problem is resolved.
\paragraph{(ii) $X \cap \mathcal {Q } \ne \emptyset $}.\\
Let us choose a point $ P \in X \cap \mathcal {Q }$. 
We need the following proposition which can be found in $\lbrack 12 \rbrack $.
\begin {Proposition}
 Let $\mathcal{Q}_{N}$ be a non-degenerate quadric in $\mathbb{P}^N$ and $H$ a hyperplane of $\mathbb{P}^N$.
Then:\\
If $H$ is tangent to $\mathcal{Q}_{N}$, then $\mathcal{Q}_{N} \cap H $ is a degenerate quadric of rank $N-1$ in  $\mathbb{P}^{N-1}$.\\ 
If $H$ is not tangent to $\mathcal{Q}_{N}$, then $\mathcal{Q}_{N} \cap H $ is a non-degenerate quadric in  $ \mathbb{P}^{N-1}$. 
\end{Proposition}
Let $H_1$ be the tangent plane to the quadric at the point $P$.  From the above proposition, we can deduce that $ \mathcal{Q} \cap H_{1}$ is a degenerate quadric in $PG(2,t^2)$ of rank equal to $2$.  In \lbrack 6\rbrack  \ p.156, table 7.2, J.W.P. Hirschfeld gave the classification of quadrics in $PG(2,q)$. Thus, we get two cases: $ \mathcal{Q} \cap H_{1}$ is either a pair of distinct lines or a point (a pair of conjuguate lines in $PG(2,q)$ which meet in $PG(2,t)$). Since the quadric $ \mathcal{Q}$  does not contain any line, we can say that $ \mathcal{Q} \cap H_{1} $ is a point, $ \mathcal{Q} \cap H_{1} =\{ P \}$ and therefore  $$  \mathcal{Q} \cap X \cap H_{1} =\{ P \}  \quad (5.2.2.1)$$ 
Let us take a line D passing through $P$ and contained in the plane $H_{1}$. 
We consider also  all the planes $(H_{i} )_{i \in I}$ passing through  the line D. \\
From the above proposition, if $H_{i} \ne H_1$, then $ \mathcal{Q} \cap H_{i} $ is a non-degenerate quadric and therefore from proposition 5.2 a non-singular curve.\\
On the other hand $X \cap H_{i} $ is a hermitian curve (singular or not). \\
If $H_{i}$ is not tangent to $X$, then $X \cap H_{i} $ is a non-singular hermitian curve. Thus, $X \cap H_{i} $ and $ \mathcal{Q} \cap H_{i} $ are two non-singular projective plane curves, so they are irreducible. From the theorem of B\'ezout  \lbrack13\rbrack \ chap.4, 2.1 p.236, we can say that $$\vert  X\cap \mathcal{Q} \cap H_{i} \vert \le 2(t+1)  \quad (5.2.2.2)$$
If $H_{i}$ is  tangent to $X$, then $X \cap H_{i} $ as a singular hermitian curve is from theorem 4.1 a set of $t+1$ lines passing through a common point; each line meeting $ X \cap \mathcal{Q} $ in at most $2$ points. Thus, we get $$\vert  X \cap \mathcal{Q} \cap H_{i} \vert \le 2(t+1) \quad (5.2.2.3)$$
We can now conclude that  $$ \mathrm{if } \ H_i \ne H_1, \mathrm{then} \qquad \vert X\cap \mathcal{Q} \cap H_{i}  \vert \le 2(t+1)  \quad (5.2.2.4 )$$
The point $P$ belong to each  $X \cap \mathcal{Q} \cap H_{i}$, and therefore from (5.2.2.4) we get $$ \vert X \cap \mathcal{Q} \cap H_{i}  -\{P\}  \vert \le 2(t+1) -1 \quad(5.2.2.5)$$
There are exactly $q+1$ planes passing through the line D, and their union generate the whole projective space $PG(3,q)$. 
Thus, we get  $$ \mathcal{Q} \cap X = \{P \} \cup \bigcup_{i=2}^{q+1}(  X \cap \mathcal{Q} \cap H_{i}  -\{P\} ) \quad (5.2.2.6)$$ 
From the relation  (5.2.2.6) we get 
$$\vert \mathcal{Q }\cap{X} \vert \le 1+ \sum_{i=2}^{q+1} {\vert X \cap \mathcal{Q} \cap H_{i}  -\{P\} \vert}_{} \qquad (5.2.2.7)$$
And from the relations (5.2.2.5) and (5.2.2.7) we deduce that $\vert \mathcal{Q }\cap X \vert \le 1+q(2(t+1)-1)$.\\
Finally we get,  $$\vert \mathcal{Q }\cap X \vert \le 2t^3+t^2+1= s_2(t).$$
It is now possible for us, to answer affirmatively to the conjecture of S\o rensen in the case $h=2$ and $t$ any prime power.
\subsection{Positive answer to S\o rensen' s conjecture in the case of quadrics}
\begin{theorem}
Let $\mathcal{Q}$ be a quadric in PG(3,q) and  X the non-degenerate hermitian surface $ X : x_{0}^{t+1}+x_{1}^{t+1}+x_{2}^{t+1}+x_{3}^{t+1}= 0$, we have  $$\# X_{ Z(\mathcal{Q})}(\mathbb{F}_{q} ) = 2( t^{3}+ t^{2}-t)+t+1\quad \mathrm{or} \qquad \# X_{ Z(\mathcal{Q})}(\mathbb{F}_{q} ) \le 2 t^{3}+ t^{2}+1   $$
  
\end{theorem}
We can now state the following results.

\section{The parameters of the functional codes defined by forms of degree 2 in $PG(3, q)$ on the hermitian surface}
\begin{theorem}
The code $C_{2}(X)$ defined on the hermitian surface $ X : x_{0}^{t+1}+x_{1}^{t+1}+x_{2}^{t+1}+x_{3}^{t+1}= 0$ is a $\lbrack n, k,d \rbrack_{q}$-code where  \\
 $$ n= (t^2+1)(t^3+1)$$ 
 $$k =10 $$
 $$d=t^5-t^3-t^2+t =t(t-1)(t^3+t^2-1)$$
\end{theorem}
\textbf{Proof}
To find $d$ we use the relations (3.2.2), (3.4.3) and the theorem 5.11.
\begin{theorem}
The minimum weight codewords correspond now to the two planes each intersecting the hermitian surface $X$ in $t+1$ lines passing trough a point and the line of intersection of the two planes intersects the hermitian surface in $t+1$ points.
\end{theorem}

\begin{lemma}\label{PP}
Let $\{\mathcal{P}_1, \mathcal{P}_2\}$ and $\{\mathcal{P}_1, \mathcal{P}_3\}$ be two distinct pairs of tangent planes to the hermitian surface $X$. Then  $(\mathcal{P}_1\cap X) \cup(\mathcal{P}_2 \cap X ) \ne  (\mathcal{P}_1\cap X) \cup(\mathcal{P}_3 \cap X ).$
\end{lemma}
\begin{lemma}   
Let $\{\mathcal{P}_1, \mathcal{P}_2\}$ and $\{\mathcal{P}_3, \mathcal{P}_4\}$ be two distinct pairs of tangent planes to the hermitian surface $X$. Then  $(\mathcal{P}_1\cap X) \cup(\mathcal{P}_2 \cap X ) \ne  (\mathcal{P}_3\cap X) \cup(\mathcal{P}_4 \cap X ).$
\end{lemma}
The proof of lemma 6.3 and lemma 6.4 are obvious. It is done by a simple calculation of intersection. \\ 

We come now to the computation of the codewords of minimum weight.
\begin{theorem}
The number of codewords of minimum weight is : $$ (t^2-1)\lbrack \frac{1}{2} (t^5+t^3+t^2+1) t^5 \rbrack $$ 
\end{theorem}
\textbf{Proof}\\
We know that a codeword of minimum weight is given by a quadric  $ \mathcal{Q}=\mathcal{P}_1\cup \mathcal{P}_2 $ where : \\
\textendash \ $\mathcal{P}_1$ is a tangent plane to the hermitian surface $X$ at a point $P_1$,\\ 
\textendash \ $\mathcal{P}_2$ is a tangent plane to the hermitian surface $X$ at a point $P_2$,\\
Ñand $P_2$  doesn't  belong in $\mathcal{P}_1$, (i.e. $ \Delta= \mathcal{P}_1 \cap \mathcal{P}_2$ meets the hermitian surface in $t+1$ points).\\
We also know from \lbrack2\rbrack, theorem 7.3 p.1172,  that there are exactly  $\# X(\mathbb{F}_q)=t^5+t^3+t^2+1$ tangent planes $ \mathcal{P}$ to the surface $X$.  Thus $\mathcal{P}_1\cap X$ as a set of $t+1$ lines meeting in a common point, has $t^3+t^2+1$ points. \\
Given a tangent plane $\mathcal{P}_1$, if we wish that the quadric $\mathcal{Q}=\mathcal{P}_1\cup \mathcal{P}_2$ gives a codeword of minimum weight, we get exactly $\# X (\mathbb{F}_q)-(t^3+t^2+1)=t^5 $ possibilities to choose the point  $P_2$. Thus, we get  $(t^5+t^3+t^2+1)t^5$  pairs  $(\mathcal{P}_1,\  \mathcal{P}_2)$ of tangent planes to the hermitian surface $X$ giving a codeword of minimum weight. Since the two pairs $(\mathcal{P}_1,\  \mathcal{P}_2)$, and  $(\mathcal{P}_2,\  \mathcal{P}_1)$ are the same (they give the same quadric), we get from these tangent planes $\frac{1}{2}(t^5+t^3+t^2+1)t^5$ quadrics giving codewords of minimum weight. From lemmas 6.3 and 6.4, the codewords obtained from these quadrics are all distinct. \\
Given a quadric $\mathcal{Q}$, from which we get a codeword of minimum weight, the quadric $\mathcal{Q}^{\prime}=\lambda \mathcal{Q}$ with $\lambda \in \mathbb{F}_{q}^{\star}$, gives a codeword of minimum weight. Thus, we get finally  $(t^2-1)\lbrack \frac{1}{2} (t^5+t^3+t^2+1) t^5 \rbrack $ codewords of minimum weight.  \\\\     
Observe that, in the case $t=2$ there is exactly $3\times \frac{1}{2}\times (45\times 32)=2160$ codewords of minimum weight. Thus, we recover the result of Paul P. Spurr \lbrack 2\rbrack \  p.120 and the one of S\o rensen \lbrack15\rbrack \ p.9.
\begin{theorem}
The second weight of the code $C_{2}(X)$ is  $t^5-t^3$.    
\end{theorem}
\textbf{Proof}
In fact when $\mathcal{Q}$ is not a pair of two tangent planes which intersection line meets $X$ in $t+1$ points, we get $\# X_{ Z(\mathcal{Q})}(\mathbb{F}_{q} ) \le s_{2}(t)=2 t^{3}+ t^{2}+1$.\\
And this upper bound is reached when $\mathcal{Q}$ is either a pair of two tangent planes meeting in a line contained in $X$, or a hyperbolic quadric containing three skew lines of $X$, or a pair of two planes, one tangent to $X$ and the second not tangent to $X$, with the line of intersection of the two planes meting $X$ at a single point.\\\\
Note that, in the case $t=2$ we recover the second weight of the code  $C_{2}(X)$ over $\mathbb{F}_4$ in agreement with  \lbrack 3\rbrack  \ which is $2^5-2^3=24$.\\\\ 
\textbf{Acknowledgments}\\\\
The author would like to thank Mr F. Rodier. He is very grateful to him. His remarks and patience encouraged me to work on the problem. \\\\
\textbf{References}\\\\\
\lbrack 1\rbrack   \ Y. Aubry,  Reed-Muller codes associated to projective algebraic varieties. In ``Algebraic Geomertry and Coding Theory ". (Luminy, France, June 17-21, 1991). Lecture Notes in Math., Vol. 1518 pp.4-17, Springer-Verlag, Berlin, 1992.\\\\
\lbrack 2\rbrack  \ R. C Bose and I. M Chakravarti, Hermitian varieties in finite projective space $PG(N,q)$. Canadian J.of Math.18 (1966), pp1161-1182.\\\\
\lbrack 3\rbrack  \ I. M. Chakravarti, The generalized Goppa codes and related discrete designs from hermitian surfaces in $PG(3, s^2)$. Lecture Notes in computer Sci 311. (1986), pp 116-124.\\\\
\lbrack 4\rbrack \  I. M. Chakravarti, Families of codes with few distinct weights from singular and nonsingular varieties and quadrics in projective geometries and Hadamard difference sets and designs associated with two-weight codes, in ``Coding Theory and Design Theory ", Part I, IMA Vol. Math. Appl. 20, Springer, New York,  1990, pp35-50. \\\\
\lbrack 5\rbrack \ R. Hartshorne, Algebraic geometry, Graduate texts in mathematics 52, Springer-Verlag, 1977. \\\\
\lbrack 6\rbrack  \ J. W. P. Hirschfeld, Projective Geometries Over Finite Fields (Second Edition) Clarendon  Press. Oxford 1998.\\\\
\lbrack 7\rbrack  \ J. W. P. Hirschfeld, Finite projective spaces of three dimensions, Clarendon press. Oxford 1985. \\\\
\lbrack 8\rbrack  \ J. W. P. Hirschfeld,  M. Tsfasman, S. Vladut, The weight hierarchy of higher dimensional hermitian codes, IEEE Trans. Inform. Theory 40 (1) 1994, pp.275-278.\\\\
\lbrack 9\rbrack \ G. Lachaud, The parameters of projective Reed-Muller codes, Discrete Math, 81 (2) (1990) 217-221.\\\\
\lbrack 10\rbrack  \ G. Lachaud, Number of points of plane sections and linear codes defined on algebraic varieties;  in " Arithmetic, Geometry, and Coding Theory " . (Luminy, France, 1993), Walter de Gruyter, Berlin-New York , 1996, pp 77-104. \\\\
\lbrack 11\rbrack \  Y. I. Manin, S. Vladut, Linear codes and modular curves, Itogi Nauki i Tekhniki 25 (1984) 209-257 (English translation J. Soviet Math. 30 (1985) 2611-2643). \\\\
\lbrack 12\rbrack \ E. J. F. Primrose,  Quadrics in finite geometries, Proc. Camb. Phil. Soc., 47 (1951), 299-304.\\\\
\lbrack 13\rbrack \  I. R. Shafarevich, Basic algebraic geometry 1, Springer-Verlag, 1994.\\\\
\lbrack 14\rbrack \ A. B. S\o rensen, Projective Reed-Muller codes, IEEE Trans. Inf. Theory Vol. 37  (6) (1991) 1567-1576.\\\\
\lbrack 15\rbrack  \ A. B. S\o rensen, Rational points on hypersurfaces, Reed-Muller codes and algebraic-geometric codes. Ph. D. Thesis, Aarhus, Denmark, 1991.

 \end{document}